\documentclass[11pt]{article}
\usepackage{amsfonts}
\usepackage[mathscr]{eucal}
 \usepackage{amsmath}

\newtheorem{theorem}{Theorem}[section]

\newtheorem{proposition}[theorem]{Proposition}
\newtheorem{definition}[theorem]{Definition}
\newtheorem{remark}[theorem]{Remark}
\newtheorem{corollary}[theorem]{Corollary}

\newcommand{\al}{{\alpha}}
\newcommand{\om}{{\omega}}
\newcommand{\Om}{{\Omega}}

\newcommand{\si}{{\sigma}}

\newcommand{\la}{\lambda}

\newcommand{\eps}{\varepsilon}
\newcommand{\ga}{{\gamma}}

\newcommand{\Ac}{\mathfrak{A}}

\newcommand{\sD}{\mathscr{D}}

\newcommand{\sW}{\mathscr{W}}

\newcommand{\R}{{\mathbb  R}}

\newcommand{\cA}{{\cal A}}
\newcommand{\cB}{{\cal B}}

\newcommand{\cO}{{\cal O}}
\newcommand{\cF}{{\cal F}}

\newcommand{\cD}{{\cal D}}
\newcommand{\cE}{{\cal E}}
\newcommand{\cH}{{\cal H}}

\newcommand{\cK}{{\cal K}}

\newcommand{\di}{{\rm div\, }}

\newenvironment{declaration}[1]{\trivlist
\item[\hskip \labelsep{\bf #1 }]\ignorespaces}{\endtrivlist}
\newenvironment{proofof}[1]{\begin{declaration}{#1}}{\hfill
$\square$ \end{declaration}}
\newenvironment{proof}{\begin{proofof}{Proof.}}{\end{proofof}}

\begin{document}
\title{A global attractor  for   a  fluid--plate interaction model
accounting only for longitudinal deformations of the plate }
\author{Igor Chueshov\footnote{\small e-mail:
 chueshov@univer.kharkov.ua}
 \\Department of Mechanics and Mathematics, \\
 Kharkov National University, \\ Kharkov, 61077,  Ukraine\\  } \maketitle
\begin{abstract}
We study asymptotic dynamics of a  coupled system consisting of linearized 3D
Navier--Stokes equations in a bounded domain  and the classical (nonlinear)
elastic plate equation for in-plane motions
 on a flexible flat  part  of the boundary.
The main peculiarity of the model is the assumption that
 the transversal displacements of the plate are negligible
relative to in-plane displacements.
This kind of models arises in the study of blood flows in large arteries.
Our main result states the existence of a compact global attractor
of finite dimension. We also show that the corresponding linearized system
generates exponentially stable $C_0$-semigroup.
We do not assume any kind of mechanical damping in the plate component.
Thus our results means that dissipation of the energy in the fluid due to viscosity
is sufficient to stabilize the system.
\end{abstract}

\section{Introduction}
We consider a coupled (hybrid) system which describes
interaction a homogeneous viscous incompressible fluid
which  occupies a domain $\cO$ bounded by
the (solid) walls of the container $S$ and a horizontal boundary $\Om$
on which a thin (nonlinear) elastic plate is placed.
The motion of the fluid is described
by  linearized 3D Navier--Stokes equations.
To describe  deformations of the plate we involve a general (full)
Kirchhoff--Karman type model (see, e.g., \cite{lagnese,lag-2,lag-lions}
and the references therein)  with additional hypotheses
that the transversal displacements of the plate are negligible
relative to in-plane displacements.
Thus we only consider longitudinal deformations of the plate
and take account of tangential shear forces
which fluid exerts on the plate.
This kind of models arises in the study of the problem
of blood flows in large arteries
(see, e.g., \cite{pedley80} and also \cite{Ggob-jmfm08}
for a more recent discussion and references).
\par
We note that the mathematical studies of the problem
of fluid--structure interaction
in the case of viscous fluids and elastic plate/bodies
have a long history.
We refer to \cite{CDEG05,Ggob-jmfm08,Ggob-aa09,Ggob-mmas09,Kop98}
and the references therein for the case of plates/membranes,
to \cite{CS06} in the case of moving elastic bodies,
and to \cite{avalos-amo07,aval-tri07,aval-tri09,BGLT07,BGLT08,DGHL03}
in the case of elastic bodies with the fixed interface;
see also  the literature cited in these references.
\smallskip\par
Our mathematical model is formulated as  follows.
\par
 Let $\cO\subset \R^3$
 be a bounded domain  with a sufficiently smooth
 boundary $\partial\cO$. We assume that $\partial\cO=\Omega\cup S$,
 where $\Om\subset\{ x=(x_1;x_2;0)\, :\,(x_1;x_2)\in\R^2\}$ and $S$
 is a smooth surface which lies under plane $x_3=0$.
 The exterior normal on $\partial\cO$ is denoted
 by $n$. We have that $n=(0;0;1)$ on $\Om$.
 We consider the following linear Navier--Stokes equations in $\cO$
for the fluid velocity field $v=v(x,t)=(v^1(x,t);v^2(x,t);v^3(x,t))$
and for the pressure $p(x,t)$:
\begin{equation}\label{fl.1}
   \partial_tv-\nu\Delta v+\nabla p=\tilde{G}\quad {\rm in\quad} \cO
   \times(0,+\infty),
\end{equation}
   \begin{equation}\label{fl.2}
   \di v=0 \quad {\rm in}\quad \cO
   \times(0,+\infty),
  \end{equation}
   \begin{equation}\label{fl.3}
   v(x,0)=v_{0}(x) \quad {\rm in}\quad
   \cO,\end{equation}
where $\nu>0$ is the dynamical viscosity and $\tilde{G}$ is a volume force.
   We supplement (\ref{fl.1})--(\ref{fl.3}) with  the (non-slip)  boundary
   conditions imposed  on the velocity field $v=v(x,t)$:
\begin{equation}\label{fl.4}
v=0 ~~ {\rm on}~S;
\quad
v\equiv(v^1;v^2;v^3)=(u_t;0)\equiv(u^1_t;u^2_t;0) ~~{\rm on} ~ \Om,
\end{equation}
where $u=u(x,t)\equiv (u^1(x,t); u^2(x,t))$ is the in-plane displacement vector
of the plate placed on $\Om$ satisfying
the following  equations  (see, e.g., \cite{lag-2} and \cite{lag-lions}):
\begin{equation}\label{pl.1}
\rho h u^1_{tt} -\left(\partial_{x_1} N_{11} +
\partial_{x_2} N_{12}\right)+ F^1=0
~~{\rm in}~~ \Omega \times (0, \infty),
\end{equation}
and
\begin{equation}\label{pl.2}
\rho h u^2_{tt} -\left(\partial_{x_1} N_{21} +
\partial_{x_2} N_{22}\right)+ F^2=0
~~{\rm in}~~ \Omega \times (0, \infty),
\end{equation}
with
\[
N_{11}=D\left(u^1_{x_1}+\mu u^2_{x_2}\right),\quad N_{22}=D\left(u^2_{x_2}+\mu u^1_{x_1}\right)
\]
and
\[
N_{12}=N_{21}=\frac{D}2(1-\mu)\left(u^1_{x_2}+ u^2_{x_1}\right),
\]
where $D=E h/(1-\mu^2)$, $E$ is Young's modulus, $0<\mu<1/2$
is Poisson's ratio, $h$ is the thickness of the plate,
$\rho$ is the mass density.
The external (in-plane) force $(F^1;F^2)$  in \eqref{pl.1} and  \eqref{pl.2}
consists of two parts,
\[
F^i= f^i(u^1,u^2) + T_i(v),\quad i=1,2,
\]
where $(f^1(u^1,u^2);f^2(u^1,u^2))$ is a nonlinear feedback force
represented by some potential $\Phi$ (which we specify below):
\[
f^i(u^1,u^2)=\frac{\partial \Phi(u^1,u^2)}{\partial u^i},\quad i=1,2,
\]
and $(T_1(v);T_2(v))$  is the viscous shear stress exerted by the fluid on the plate,
$T_i(v)=\left( (Tn) \big|_\Om, e_i\right)_{\R^3}$.  Here $T=\{T_{ij}\}_{i,j=1}^3$
is the stress tensor of the fluid,
\[
T_{ij}\equiv T_{ij}(v)=\nu\left(v^i_{x_j}+v^j_{x_i}\right)-p\delta_{ij},
\quad i,j=1,2,3,
\]
$e_1=(1;0;0)$,  $e_2=(0;1;0)$ are unit tangential vectors
on $\Om\subset\partial\cO$
and $n=(0;0;1)$
is the outer normal vector to $\partial\cO$ on $\Om$.
A simple calculation shows that
\[
T_{i}(v)\equiv T_{i3}(v)=\nu\left(v^i_{x_3}+v^3_{x_i}\right) =\nu v^i_{x_3},
\quad i=1,2,
\]
(in the last equality we use the fact that
$v^3(x_1;x_2;0)=0$ for $(x_1;x_2)\in\Om$ due to the second relation in \eqref{fl.4}
and hence $v^3_{x_i}=0$ on $\Om$, $i=1,2$).
\par
Thus we arrive at the following equations for  the in-plane displacement
$u=(u^1; u^2)$  of the plate
(below for some notational simplifications we assume that $\rho h=1$
and $D(1-\mu)/2=1$):
\begin{equation}\label{pl-1m}
                        u^i_{tt} -\Delta u^i-\la \partial_{x_i}\left[ {\rm div}\, u\right]
+\nu    v^i_{x_3}|_{x_3=0}  +
                        f^i( u)=0, \quad i=1,2,
\end{equation}
where $\la=(1+\mu)(1-\mu)^{-1}$ is a nonnegative parameter.
For  the  displacement $u=(u^1; u^2)$
we impose  the  clamped boundary conditions  on
$\Gamma=\partial\Om$:
\begin{eqnarray}\label{pl-2m}
u^i =0~~{\rm on}~~  \Gamma, ~~ i=1,2.
\end{eqnarray}
Our main point of interest is wellposedness and long-time dynamics of
solutions to the coupled problem in
\eqref{fl.1}--\eqref{fl.4}, \eqref{pl-1m}, and \eqref{pl-2m}
for the velocity $v$ and the displacement $u=(u^1; u^2) $ with the initial data
\begin{equation}\label{fp-1}
v|_{t=0}=v_0, ~~ u|_{t=0}=u_0,~~ u_t|_{t=0}=u_1.
\end{equation}
This problem in the case when $\tilde G\equiv 0$, $\la=0$ and $f^i(u)\equiv 0$
was considered in \cite{Ggob-mmas09} (see also \cite{Ggob-jmfm08,Ggob-aa09})
with an additional strong (Kelvin-Voight type)
damping force applied to the interior of the plate. These papers deals with
the existence and asymptotic stability of the corresponding semigroup.
In contrast with \cite{Ggob-jmfm08,Ggob-aa09,Ggob-mmas09}
we do not assume the presence of mechanical damping terms in the plate component
of the system and consider nonlinearly forced model.
\par
Our main result (see Theorem~\ref{th:main}) states that under some natural conditions concerning
the potential $\Phi(u)$ of feedback forces system
\eqref{fl.1}--\eqref{fl.4} and \eqref{pl-1m}--\eqref{fp-1}
possesses a compact global attractor of finite fractal dimension
which also has some additional regularity properties.
Our considerations involves   recently
developed approach based on quasi-stability properties
and stabilizability estimates
(see, e.g., \cite{cl-jdde,cl-mem,cl-book} for the second order
in time evolution equations and also
\cite{BC-08,cl-amo08} for thermoelastic problems).
\par
We also improve substantially  the stability results
presented in \cite{Ggob-jmfm08,Ggob-aa09,Ggob-mmas09}.
Namely, as a consequence of dissipativity property
of the nonlinear model we prove (see Corollary~\ref{co:as-stab})
that linear part of \eqref{fl.1}--\eqref{fl.4} and \eqref{pl-1m}--\eqref{fp-1}
generates uniformly exponentially stable
$C_0$-semigroup of contractions
 {\em without} any dissipation mechanisms (except the fluid viscosity)
in the system.
\par
We note that a linear model related to
\eqref{fl.1}--\eqref{fl.4} and \eqref{pl-1m}--\eqref{fp-1}
was considered in \cite[Section 3.15]{lt-book} from the point of view
of boundary control.
In \cite{lt-book} the authors deal with the so-called
linearized model
of well/reservour coupling with monophasic flow.
This is a scalar linear model represented by
the diffusion (heat) equation in $\cO$ coupled with
the wave equation  on $\Om$ with an interface
condition like in \eqref{fl.4}.
The scalar structure of the model makes it possible to prove
 \cite[Proposition 3.15.5]{lt-book}
that this model generates contractive exponentially stable analytic semigroup.
We do not know whether the result on analyticity remains true for the hydrodynamical
situation we consider.
We also note that in fact  the long time dynamics  in problem \eqref{pl-1m}
with a given {\em stationary smooth} velocity field $v$ was studied in
\cite{cl-mt} and \cite[Chapter 7]{cl-mem}.
\par
The paper is organized as follows. In the next Section~\ref{sec:pre}
we discuss an auxiliary Stokes problem and rewrite problem
\eqref{fl.1}--\eqref{fl.4} and \eqref{pl-1m}--\eqref{fp-1}
as a single first order equation in some extended space $\cH$.
Then we prove that the linear version of the problem
generates strongly continuous  contractive semigroup in $\cH$.
In Section~\ref{sec:wp} we prove the existence of strong and generalized
solutions to the original nonlinear problem
and establish energy balance identity.
For this we use ideas presented in \cite{pazy}.
 Section~\ref{sec:main} contains our main results on
 long-time dynamics  of system
\eqref{fl.1}--\eqref{fl.4} and \eqref{pl-1m}--\eqref{fp-1}.
Our considerations here are based mainly on the idea of quasi-stability
(see, e.g., \cite{cl-mem} and \cite[Section~7.9]{cl-book}).

\section{Preliminaries}\label{sec:pre}
In this section we first provide with some results concerning
to  Stokes problem and then consider the abstract form of the system.
Below we denote by $H^s(D)$ the Sobolev space of the order $s$
on the set $D$, by $H^s_0(D)$ the closure of $C_0^\infty(D)$ in  $H^s(D)$,
and by  $H^s(D)/\R$ the factor-space with the naturally induced norm.
\subsection{Stokes problem}

In further considerations we need some regularity  properties
of the terms responsible for fluid--plate interaction.
For this we consider
 the following Stokes problem
\begin{eqnarray}\label{stokes}
 &&  -\nu\Delta v+\nabla p= g, \quad
   \di v=0 \quad {\rm in}\quad \cO;
\\
&& v=0 ~~ {\rm on}~S;
\quad
v=(\psi;0) ~~{\rm on} ~ \Om, \nonumber
\end{eqnarray}
where $g\in [L^2(\cO)]^3$ and $\psi=(\psi^1;\psi^2)\in [L^2(\Om)]^2$ are given.
This type of boundary value problems for the Stokes equation
was studied by many authors  (see, e.g., \cite{lad-NSbook} and \cite{temam-NS}
and the references therein). We collect some properties of solutions
to \eqref{stokes} in the following assertion.
\begin{proposition}\label{pr:stokes}
Let $g\in [L^2(\cO)]^3$ and $\psi\in [H^1_0(\Om)]^2$.
Then \begin{itemize}
\item [{ \bf (1)}]
Problem \eqref{stokes} has a unique solution
$\{v;p\}\in [H^{3/2}(\cO)]^3\times[ H^{1/2}(\cO)/\R]$ such that
\begin{equation}\label{stokes-bnd1}
\|v\|_{[H^{1+s}(\cO)]^3}+\|p\|_{H^{s}(\cO)/\R}
\le c_0\left\{\|g\|_{[H^{-1+s}(\cO)]^3}+\|\psi\|_{[H^{s+1/2}(\Om)]^2} \right\}
\end{equation}
for every $0\le s\le 1/2$. Moreover,
  \item [{ \bf (2)}] We have that $\ga_\Om v\in [H^{-1/2}(\Om)]^2$
for the
 trace operator
$\ga_\Om$  defined (on smooth functions) by the formula
\begin{equation}\label{gamma-om}
\ga_\Om v=\nu \left(v^1_{x_3};v^2_{x_3}\right)
~~\mbox{for}~~v=(v^1;v^2;v^3)\in  [H^2(\cO)]^3.
\end{equation}
  \item [{ \bf (3)}]If $g=0$ we also have  that
  \begin{equation}\label{stokes-bnd2}
\|v\|_{[L_2(\cO)]^3}\le c_0 \|\psi\|_{[L_2(\Om)]^2},
\end{equation}
thus we can define a bounded operator
$N_0 :\left[L_2(\Om)\right]^2\mapsto [L^2(\cO)]^3$ by the formula
\begin{equation}\label{fl.n0}
N_0\psi=w ~~\mbox{iff}~~\left\{
\begin{array}{l}
 -\nu\Delta w+\nabla p=0, \quad
   \di w=0 \quad {\rm in}\quad \cO;
\\
 w=0 ~~ {\rm on}~S;
\quad
w=(\psi;0) ~~{\rm on} ~ \Om,
\end{array}\right.
\end{equation}
for $\psi=(\psi^1;\psi^2)\in \left[L_2(\Om)\right]^2$
($N_0\psi$ solves \eqref{stokes} with $g\equiv 0$).
\end{itemize}
\end{proposition}
\begin{proof} {\bf 1.}
The existence and uniqueness of solutions along with the bound in
\eqref{stokes-bnd1} follow from Proposition 2.3 and Remark 2.6
on  Sobolev norm's interpolation in \cite[Chapter~1]{temam-NS}.
\par
{\bf 2.} To prove Statement 2 we use the same idea as in \cite{avalos-amo07}
which involves the boundary properties
of harmonic functions
(see \cite{kellogg}).
\par
We can represent $v$ in the form $v=\hat v+v^*$,
where $\hat v$ solves \eqref{stokes} with $\psi\equiv 0$
and $v^*$  satisfies \eqref{stokes}  with $g\equiv 0$.
By Proposition 2.3~\cite[Chapter~1]{temam-NS}.
we have that $\hat v\in [H^{2}(\cO)]^3$
and thus by the standard trace theorem there exists
$\partial_n\hat v|_{\partial\cO}\in [H^{1/2}(\partial\cO)]^3$.
Consequently $\ga_\Om\hat v\in [H^{1/2-\delta}(\Om)]^2
\subset [H^{-1/2}(\Om)]^2$.
\par
Thus we need to establish Statement 2 in the case $g\equiv 0$ only.
In this case the pressure $p$ is a harmonic function in $\cO$
which belongs $L_2(\Om)$. Thus
by Theorem~3.8.1~\cite{kellogg}
we can assign a meaning to $p|_\Om$ in $H^{-1/2}(\partial\cO)$.
Now using the Gauss-Ostrogradskii formula one can see
that $\nabla p\in \left([H^1(\cO)]^3\right)'$.
Therefore from \eqref{stokes} with $g=0$
 we have $\Delta v\in \left([H^1(\cO)]^3\right)'$.
 Using the Green formula
 \[
\int_{\partial\cO} \partial_n v ^i\psi dS=
\int_{\cO} \Delta v ^i\psi dx- \int_{\partial\cO}\nabla v ^i\nabla \psi dx,
\quad \forall \psi\in H^1(\cO),
\]
for every velocity component $v^i$,
we conclude that $\ga_\Om v\in [H^{-1/2}(\Om)]^2$.
\par
{\bf 3.}
To prove the  third statement we use the representation
of the velocity field with hydrodynamical potentials
(see \cite{lad-NSbook}).
Indeed (see \cite[Chapter 3, Section 2]{lad-NSbook}),
the velocity field $v=(v^1;v^2;v^3)$ solving \eqref{stokes}
with $g\equiv 0$  admits the representation
\begin{equation}\label{pot-1}
v^i(x) =\sum_{j=1}^3 \int_{\partial\cO} K_{ij}(x,y) \phi^j(y) dS_y,~~ i=1,2,3, ~ x\in\cO,
\end{equation}
where  $\phi=(\phi^1;\phi^2;\phi^3)$ satisfies the integral equation on $\partial\cO$:
\begin{equation}\label{pot-2}
\frac12 \phi^i(x)+\sum_{j=1}^3 \int_{\partial\cO} K_{ij}(x,y) \phi^j(y) dS_y
=\tilde\psi^i(x),~~ i=1,2,3,~ x\in \partial\cO.
\end{equation}
Here $\tilde\psi$ is the extension of $(\psi;0)$ on $\partial\cO$ by zero.
The kernels $K_{ij}(x,y)$ has the form
\[
K_{ij}(x,y)=-\frac{3}{4\pi}\frac{(x_i-y_i)(x_j-y_j)}{|x-y|^2}\cdot
\frac{(x-y, n_y)}{|x-y|^3},~~ x\in\R^3,~ y\in \partial\cO,
\]
where $n_y$ is the outer normal vector at $y$.
By the theory of integral equations with weak singularities
(see, e.g., the references in \cite{lad-NSbook}) we have that
the function $\phi$ is uniquely defined by $\psi$
from \eqref{pot-2}
and the mapping $\psi\mapsto\phi$ is continuous
from $[L_2(\Om)]^2$ into $[L_2(\partial\cO)]^3$.
One can also show (some details can be found in
\cite[Chapter 3, Section 2]{lad-NSbook}) that
\[
\sup_{x\in\cO}\int_{\partial\cO}| K_{ij}(x,y)|dS_y +
\sup_{y\in\partial\cO}\int_{\cO}| K_{ij}(x,y)|dx\le C,~~ i,j=1,2,3.
\]
Therefore we can apply the Schur test (see, e.g., Theorem~5.2 in \cite{halmos-sanders})
to conclude that the mapping $\phi\mapsto v$ given by \eqref{pot-1}
is continuous from $[L_2(\partial\cO)]^3$
into $[L_2(\cO)]^3$.
This implies \eqref{stokes-bnd2}.
\end{proof}
We do not pretend that the boundary trace regularity stated
in Proposition~\ref{pr:stokes} are optimal. They are sufficient for our purposes
and therefore we do not pursue the optimality issues.

\subsection{Abstract representation of the problem}

We introduce the following spaces
\[
X=\left\{ v\in [L_2(\cO)]^3\, :\; {\rm div}\, v=0;\;
\gamma_n v\equiv (v,n)=0~\mbox{on}~ \partial\cO\right\};
\]
and
\[
V=\left\{ v\in [H^1(\cO)]^3\, :\; {\rm div}\, v=0;\;
v=0~\mbox{on}~ S,\;  \gamma_n v=0~\mbox{on}~\Om \right\}.
\]
We equip   $X$ with $L_2$-type norm $\|\cdot\|_\cO$
and denote by $(\cdot,\cdot)_\cO$ the corresponding inner product.
The space $V$ is endowed  with the norm  $\|\cdot\|_V= \|\nabla\cdot\|_\cO$
\par
Let $P_L$ be the Leray projector which maps
 the space $[L_2(\cO)]^3$ onto $X$. With this projector we rewrite
 \eqref{fl.1}--\eqref{fl.4} in $X$ as follows

\begin{equation}\label{fl-abstr}
   \partial_tv+A_0(v-N_0u_t) =G,~~ t>0,\quad v|_{t=0}=v_0,
\end{equation}
where $G=P_L\tilde{G}\in X$, $N_0$ is defined by \eqref{fl.n0}, and $A_0= -\nu P_L\Delta$
is a positive self-adjoint operator
with the domain
\[
\cD(A_0)=\left\{ v\in [H^2(\cO)]^3\, :\; {\rm div}\, v=0;\;
v=0~\mbox{on}~ \partial\cO \right\}.
\]
\par
For the plate component  we use the spaces
\[
Y= \left[L_2(\Om)\right]^2\equiv L_2(\Om)\times L_2(\Om)
~~\mbox{and}~~ W=\left[H^1_0(\Om)\right]^2.
\]
We denote by $\|\cdot\|_\Om$
and  $(\cdot,\cdot)_\Om$ the norm and the inner product in $Y$.
\par
In the space $Y$ for $\la>0$ we introduce the operator
\begin{equation}\label{opA-pl}
A=- \left[
\begin{array}{cc}
(1+\la)  \partial_{x_1}^2+ \partial_{x_2}^2  & \la
\partial_{x_1x_2}^2 \\
  \\
 \la \partial_{x_1x_2}^2& \partial_{x_1}^2+ (  1+\la)\partial_{x_2}^2
\\
\end{array} \right]
\end{equation}
with the domain  $\cD(A)= [(H^2\cap H^1_0)(\Om)]^2\subset W$.
One can see that $A$ is  a positive operator in $Y$ generated by the form
\[
a(u, \hat u)=\sum_{i=1}^2\int_\Om\nabla u^i \nabla \hat u^i d\Om+
\la \int_\Om\di u\cdot \di \hat u\, d\Om \equiv (A^{1/2} u, A^{1/2}\hat u)_\Om,
\]
where $u=(u^1;u^2)$ and $u=(\hat u^1;\hat u^2)$ are from $W$.
With this operator $A$ problem \eqref{pl-1m} and \eqref{pl-2m}
can be written in the space $Y$ as follows,
\begin{equation}\label{pl-abstr}
                        u_{tt} +A u +\ga_\Om v+f(u)=0,~~ t>0, \quad u|_{t=0}=u_0,~~
                        u_t|_{t=0}=u_1.
\end{equation}
where $u=(u^1;u^2)$, the operator $A$ is defined in  \eqref{opA-pl},
$\ga_\Om$ is given by \eqref{gamma-om},
and $f(u)=(f^1(u);f^2(u))$.
\smallskip\par
Now we  consider the phase space $\cH=X\times W\times Y$
with the inner product
\[
(U, U^*)_\cH= (v,v^*)_\cO +(A^{1/2}u_0,A^{1/2}u^*_0)_\Om+(u_1,u^*_1)_\Om,
\]
where  $U=(v;u_0;u_1)$ and $U^*=(v^*;u^*_0;u^*_1)$ are elements from $\cH$.
\par
We
rewrite  problem
\eqref{fl-abstr} and \eqref{pl-abstr}
as a first order equation for the phase variable
$U=(v;u;u_t)\in \cH$ of the form
\begin{equation}\label{1st-ord}
\frac{dU}{dt}+\cA U+ \cF(U)=0,~~ t>0,\quad U|_{t=0}= U_0,
\end{equation}
where $\cF(U)=(0;0;f(u))$ and
\begin{equation}\label{opA-1}
\cA=\left[\begin{array}{ccc}
      -\nu P_L\Delta  & 0 & 0 \\
      0 & 0 & -I \\
      \ga_\Om & A & 0
    \end{array}\right]\equiv \left[\begin{array}{ccc}
      A_0  & 0 & -A_0N_0 \\
      0 & 0 & -I \\
      \ga_\Om & A & 0
    \end{array}\right]
\end{equation}
with the domain
\begin{equation}\label{opA-2}
\sD(\cA)=\left\{U=\left. \left[\begin{array}{c}
      v \\
      u_0 \\
      u_1
    \end{array}\right]\in \cH\; \right|  \begin{array}{l}
                                 A_0(v-N_0u_1)\equiv -\nu P_L\Delta v  \in X \\
                                 u_1\in W \\
                                 Au_0+\ga_\Om v\in Y,~ v|_\Om=(u_1;0)
                               \end{array}
     \right\}
\end{equation}
To solve \eqref{1st-ord} we use methods from \cite{pazy}.
For this we first prove the following assertion.

\begin{proposition}\label{pr:cA-generator}
The operator $\cA$ is a maximal accretive operator in $\cH$.
Moreover, $R(\cA)=\cH$ and thus $\cA$ is invertible and
by the Lumer--Phillips theorem (see \cite[p. 14]{pazy}) generates $C_0$-semigroup
of contractions in $\cH$.
\end{proposition}
\begin{proof}
One can see that
\[
(\cA U, U^*)_\cH= \nu(\nabla v,\nabla v^*)_\cO -(Au_0,u^*_1)_\Om+(Au^*_0,u_1)_\Om,
\]
where  $U=(v;u_0;u_1)$ and $U^*=(v^*;u^*_0;u^*_1)$
are elements from $\sD(\cA)$. This implies that
$(\cA U, U)_\cH= \nu\|\nabla v\|^2_\cO\ge 0$ and thus
$\cA$ is accretive.
\par
To prove maximality it suffices to show that $R(\cA)=\cH$, i.e., to solve
the equation of the form $\cA U= F$ for $F\equiv(g; h_0;h_1)\in \cH$.
 We obviously have that
$U=(v;u;-h_0)$, where $v\in V$ solves
 \begin{eqnarray}\label{st-eq1}
 &&  -\nu\Delta v+\nabla p= g, \quad
   \di v=0 \quad {\rm in}\quad \cO;
\\
&& v=0 ~~ {\rm on}~S;
\quad
v=(-h_0;0) ~~{\rm on} ~ \Om, \nonumber
\end{eqnarray}
with $h_0\in W=\left[H^1_0(\Om)\right]^2$,
and $u\in W$ satisfies the equation
 \begin{eqnarray}\label{st-eq2}
 A u= - \ga_\Om v+ h_1,\quad  h_1\in Y=\left[L_2(\Om)\right]^2.
\end{eqnarray}
It follows from Proposition~\ref{pr:stokes} that there exists $v\in [H^{3/2}(\cO)]^3$
satisfying \eqref{st-eq1} such that
 $\ga_\Om v\in [H^{-1/2}(\Om)]^2\subset [\cD(A^{1/4})]'$. Now we can  solve equation
\eqref{st-eq2} with respect to $u$. Thus $R(\cA)=\cH$.
\end{proof}
\begin{remark}\label{re:A-dom}{\rm
The argument above shows that
\[
\sD(\cA)\subset \left\{U=\left. \left[\begin{array}{c}
      v \\
      u_0 \\
      u_1
    \end{array}\right]\in \cH\; \right|  \begin{array}{l}
                                 v\in [H^{3/2}(\cO)]^3,\ga_\Om v\in [H^{-1/2}(\Om)]^2,
                                 \\v|_S=0,
                                    ~ v|_\Om =(u_1;0),\\
                                 u_0\in  [(H^{3/2}\cap H^1_0)(\Om)]^2,~
                                  u_1\in  [H^1_0(\Om)]^2
                               \end{array}
     \right\}\, .
\]
}
\end{remark}

\section{Well-posedness}\label{sec:wp}
We assume that  the plate force potential $\Phi(u)$ possesses the properties
\begin{itemize}
  \item $\Phi(u)\in C^2(\R^2)$ is nonnegative polynomially bounded function, i.e.,
\begin{equation}\label{phi-1}
\left|\frac{\partial\Phi(u)}{\partial u^i \partial u^j}\right|\le
C\left( 1+|u|^{p}\right), ~~ i,j=1,2,~~u= (u^1, u^2)\in \R^2,
\end{equation}
for some $C, p\ge 0$.
  \item The following dissipativity condition holds:
  for any $\delta>0$ there exists $c_1(\delta)>0$ and $c_2(\delta)\ge 0$
  such that
\begin{equation}\label{phi-2}
\sum_ {i=1,2}u^if^i(u)- c_1(\delta)\Phi(u) +
\delta |u|^2\ge -c_2(\delta) ~~\mbox{with} ~~ f^i(u)=\frac{\partial\Phi(u)}{\partial u^i}.
\end{equation}
\end{itemize}
As  examples of a such potential $\Phi(u)$ we can consider
\begin{equation}\label{phi-exm}
\Phi(u)=\psi_0(|u^1|^2+ |u^2|^2)~~\mbox{or}~~
\Phi(u)=\psi_1(u^1) +\psi_2(u^2),
\end{equation}
where $\psi_i(s)$ are nonnegative functions from $C^2(\R)$ such that
\begin{enumerate}
  \item[(a)] $|\psi''(s)|\le C\left( 1+|u|^{q}\right)$ for  some $C, q\ge 0$;
  \item[(b)]
  $s\psi'_i(s)-c_0\psi_i(s) \ge - c_1$ for some $c_0>0$ and $c_1\ge 0$.
\end{enumerate}
 For instance, $\psi_i(s)$ can be
polynomials of even degree with
a positive leading coefficient and
and with sufficiently large free term.
\par
We also note that under the assumption in  \eqref{phi-1}
the nonlinear force
\[
f(u)=(f^1(u); f^2(u))=\left(\frac{\partial\Phi(u)}{\partial u^1};
\frac{\partial\Phi(u)}{\partial u^2}\right)
\]
 in \eqref{pl-abstr} satisfies the following
local Lipschits property
\begin{equation}\label{Lip-f}
\| f(u)-f(\hat u)\|_Y\le
C_\si\left( 1+ \| u\|^p_{1-\si/p,\Om} +\|\hat u\|^p_{1-\si/p,\Om} \right)
 \| u-\hat u\|_{\si,\Om},
\end{equation}
for every $0<\si<1$ and $u,\hat u\in W$, where $\|\cdot\|_{\si,\Om}$
denotes the norm in the space $[H^{\si}(\Om)]^2$. To see this
we note that \eqref{phi-1}
implies the estimate
\[
| f^i(u)-f^i(\hat u)|\le
C\left( 1+ | u|^p +|\hat u|^p \right)| u-\hat u|,
 ~~ u,\hat u\in \R^2.
\]
Therefore the H\"older inequality and the embedding $H^{\si}(\Om)\subset L_{2/(1-\si)}(\Om)$
for  $0<\si<1$
imply \eqref{Lip-f}.
\par
Following the standard semigroup approach (see, e.g., \cite{pazy} or \cite{show})
we give the following definition.

\begin{definition}\label{de:solution}{\rm
The function $U(t)=(v(t);u(t);u_t(t))\in C(0,T;\cH)$  such that
$U(0)=U_0=(v;u_0;u_1)$
is said to be
\begin{itemize}
  \item {\em strong} solution to problem \eqref{1st-ord} on an interval
  $[0,T]$ if (i)  $U(t)\in \cD(\cA)$ for almost all $t\in [0,T]$, (ii)
  $U_t\in L_1(0,T;\cH)$, and (iii) \eqref{1st-ord}
  is satisfied as an equality in $\cH$ for almost all $t\in [0,T]$;
  \item {\em generalized } solution to problem \eqref{1st-ord} if there exist
  a sequence of initial data $U_0^n$ and the corresponding strong solutions
$U^n(t)$ such that
\[
\lim_{n\to 0}\max_{t\in [0,T]}\| U^n(t)-U(t)\|_\cH=0.
\]
\end{itemize}
}
\end{definition}

\begin{theorem}\label{th:wp}
Let $U_0\in \cH$. Then for any interval $[0,T]$
 there exists a uni\-que generalized solution
$U(t)=(v(t); u(t); u_t(t))$
such that $v\in L_2(0,T; V)$ and
 the energy balance equality
\begin{equation}\label{energy}
\cE(v(t), u(t), u_t(t))+\nu \int_0^t\|\nabla v\|_\cO^2 d\tau=\cE(v_0, u_0, u_1)
+\int_0^t(G,  v)_\cO d\tau
\end{equation}
holds for $t>0$,
where
\begin{equation}\label{en-f}
\cE(v(t), u(t), u_t(t))=\frac12\|v(t)\|^2_\cO+E(u(t), u_t(t))
\end{equation}
with the energy $ E(u,u_t)$ of the plate given by
\[
E(u,u_t)=\frac12\left( \|u_t\|^2_\Om+\| A^{1/2}u\|_\Om^2\right)+ \int_\Om\Phi(u(x))d\Om.
\]
Moreover,
\begin{itemize}
  \item Any generalized generalized solution is also mild, i.e.
\begin{equation}\label{mild-s}
U(t)=e^{-t\cA}U_0+\int_0^t e^{-(t-\tau)\cA}\cF(U(\tau)) d\tau, ~~ t>0.
\end{equation}
  \item There exists a constant $a_{R,T}>0$ such that
  for any couple of generalized solutions
  $U(t)=(v(t); u(t); u_t(t))$ and $\hat U(t)=(\hat v(t); \hat u(t); \hat u_t(t))$
with the initial data  possessing the property  $\|U_0\|_\cH, \|\hat U_0\|_\cH\le R$
we  have
\begin{equation}\label{mild-dif}
\|U(t)-\hat U(t)\|^2_\cH+\int_0^t\|\nabla( v-\hat v)\|_\cO^2 d\tau \le a_{R,T}
\|U_0-\hat U_0\|^2_\cH ~~
\end{equation}
for every $t\in [0,T]$.
  \item The  solution $U(t)$  is strong if $U_0\in \sD(\cA)$.
\end{itemize}
\end{theorem}

\begin{proof} By Proposition~\ref{pr:cA-generator} the linear part of
\eqref{1st-ord} generates a strongly continuous semigroup.
By \eqref{Lip-f} we also have that the nonlinear part in \eqref{1st-ord}
is locally Lipschitz on $\cH$. Therefore
the existence of (local) strong and generalized solutions
follows from Theorem 6.1.6 \cite{pazy}
and from the argument provided in the proof of Theorem~6.1.4 \cite{pazy}.
The latter theorem also means that generalized solutions  are mild
on the existence interval.
\par
Now we consider strong solutions on the existence interval
  and  establish the energy relation.
A simple calculation gives
\begin{equation}\label{2.2}
\frac12\frac{d}{dt}\|v(t)\|^2_\cO+\nu \|\nabla v(t)\|^2_\cO=
\nu \int_\Om \left[\frac{\partial v^1}{\partial n} u^1_t+
\frac{\partial v^1}{\partial n} u^2_t\right] d\Om
+\int_\cO vGdx.
\end{equation}
It is also clear that for plate we have that
\begin{equation}\label{2.3}
\frac{d}{dt}E(u(t),u_t(t))+
\nu \int_\Om \left[\frac{\partial v^1}{\partial x_3} u^1_t+
\frac{\partial v^1}{\partial x_3} u^2_t\right] d\Om=0
\end{equation}
The sum of equalities (\ref{2.2}) and (\ref{2.3}) along with
the relation $n=(0;0;1)$ on $\Om$ after integration  gives the
energy equality in \eqref{energy} for
strong solutions on the existence interval.
\par
It follows from \eqref{energy} and \eqref{phi-1} that
\begin{equation}\label{2.3a}
\|U(t)\|^2_\cH\le C_R + \int_0^t\| G(\tau)\|_\cO^2 d\tau,~~
\|U_0\|_\cH\le R,
\end{equation}
on the existence interval.
This implies (see, e.g., Theorem~6.1.4 \cite{pazy})
that both strong and generalized solutions cannot blow up and therefore
they can be extended on $\R_+$.
\par
By the same argument as in the proof of the energy
equality  using relations \eqref{Lip-f} and \eqref{2.3a} we can prove that
\begin{eqnarray*}
 & &\|U(t)-\hat U(t)\|^2_\cH+2\nu\int_0^t\|\nabla( v-\hat v)\|_\cO^2 d\tau
\\
& \le &
\|U_0-\hat U_0\|^2_\cH  + c_{R,T} \int_0^t\|U(t)-\hat U(t)\|_\cH \|u_t-\hat u_t\|_\Om d\tau,
~~t\in [0,T].
\end{eqnarray*}
Using the standard trace theorem from the boundary condition
in \eqref{fl.4} we have that
$\|u_t-\hat u_t\|_\Om\le C \|\nabla( v-\hat v)\|_\cO$. Therefore applying
the Gronwall type argument we obtain \eqref{mild-dif}.
\par
The relation in \eqref{mild-dif} allows us  to make limit transition
 in energy relation \eqref{energy} from strong to generalized solutions and
to conclude the proof of Theorem~\ref{th:wp}.
\end{proof}
Theorem~\ref{th:wp} makes it possible to define a
dynamical system $\big(\cH, S_t\big)$  with the phase
space $\cH=X\times W\times Y$ and the evolution operator
$S_t$ defined by
\[
S_tU_0=U(t)=(v(t);u(t);u_t(t)),
\]
where $U(t)$ is  a generalized solution to  \eqref{1st-ord} with initial data $U_0\in \cH$.

\section{Main result}\label{sec:main}
Our main result is the following theorem.
\begin{theorem}\label{th:main}
The dynamical system $(\cH, S_t)$ generated by \eqref{1st-ord} possesses a compact global
attractor $\Ac$ of finite fractal dimension. Moreover,
for any trajectory $\{ U(t)=(v(t);u(t);u_t(t)): t\in \R \}$ from the attractor $\Ac$
we have that $U\in C(\R;\sD(\cA))\cap C^1(\R;\cH)$ and
\begin{eqnarray}\label{reg}
&&  \sup_{t\in\R}\Big\{ \|v_t(t)\|_\cO + \|P_L\Delta v(t)\|_\cO\nonumber
  \\ && {}\quad +\, \| A^{1/2}u_t(t)\|_\Om
   +\| u_{tt}(t)\|_\Om+\| Au(t)+\ga_\Om v(t)\|_\Om \Big\}
\le C_\Ac.
\end{eqnarray}
\end{theorem}
We recall  (see, e.g., \cite{BabinVishik, Chueshov,Temam}) that
a \textit{global attractor}
of a dynamical system  $\big(\cH, S_t\big)$ is a bounded closed  set $\Ac\subset \cH$
 which is  invariant (i.e., $S_t\Ac=\Ac$)
 and  uniformly  attracts all other bounded  sets:
$$
\lim_{t\to\infty } \sup\{{\rm dist}(S_tU,\Ac):\ U\in B\} = 0
\quad\mbox{for any bounded  set $B$ in $\cH$.}
$$
The {\em fractal dimension} $\dim^\cH_f \Ac$ of a  set $\Ac$
in the space $\cH$ is defined as
\[
\dim^\cH_f\cA=\limsup_{\eps\to 0}\frac{\ln N(\Ac,\eps)}{\ln (1/\eps)}\;,
\]
where $N(\Ac,\eps)$ is the minimal number of closed sets in $\cH$ of
diameter $2\eps$ needed to cover the set $\Ac$.
\par
To prove the existence of  a compact global attractor for $(\cH, S_t)$
it is sufficient (see, e.g., \cite{Ha88,Temam})
to show that the system $(\cH, S_t)$ is dissipative and asymptotically smooth.
We recall (see \cite{BabinVishik, Chueshov,Temam}) that the system  is {\em dissipative} if
 there exists a bounded absorbing set $\cB_0$ in $\cH$.
A set $\cB_0$ is said to be {\em absorbing} for  $(\cH, S_t)$
if for any bounded set $\cB\subset \cH$ there exists
time $t_\cB$ such that $S_t\cB\subset \cB_0$ for all $t\ge t_\cB$.
A system $(\cH, S_t)$ is said to be \textit{asymptotically smooth}
if for any  closed bounded forward invariant set $B\subset \cH$
there exists a compact set $\cK=\cK(B)$ which  uniformly attracts $B$:
$$
\lim_{t\to\infty} \sup\{{\rm dist}(S_tU,\cK):\ U\in B\} = 0.
$$
Below to prove dissipativity we use an appropriate Lyapunov function. As fo
 asymptotic smoothness of the system and
finite-dimensionality of the attractor, we rely on recently developed approach
based on stabilizability estimates (see \cite{cl-mem,cl-book}
and the references therein).

\subsection{Dissipativity}

\begin{proposition}\label{pr:diss}
The system $(\cH, S_t)$  is dissipative.  Moreover, there exists a
 bounded  forward invariant absorbing set.
\end{proposition}
\begin{proof}
Let $S_tU_0=U(t)=(v(t);u(t);u_t(t))$.
We consider the following Lyapunov type function
\[
\sW(v, u, u_t)= \cE(v, u, u_t)+\eta \left[ (u,u_t)_\Om +(v, N_0 u)_\cO\right],
\]
where the energy $\cE$ is defined by \eqref{en-f}
and the operator $N_0$ id given by \eqref{fl.n0}.
The parameter $\eta$ will be chosen later.
\par
 In the further calculations
we deal with strong solutions.
\par
The trace theorem and the boundary condition on $\Om$ given
in \eqref{fl.4} imply that $\|u_t\|_\Om^2\le C\|\nabla v\|_\cO^2$ and
therefore it follows from energy relation \eqref{energy}
that
\begin{equation}\label{energy-d1}
\frac{d}{dt}\cE(v(t), u(t), u_t(t))+c_0( \|u_t\|_\Om^2+\|\nabla v\|_\cO^2)
\le c_1  \|G\|_\cO^2
\end{equation}
with some positive constants $c_i$.  Using \eqref{fl-abstr} and \eqref{pl-abstr}
  we also have that
\begin{eqnarray*}
& &\frac{d}{dt} \left[ (u,u_t)_\Om +(v, N_0 u)_\cO\right] \\
&= & \|u_t\|^2_\Om  + (u,u_{tt})_\Om +(v_t, N_0 u)_\cO+(v, N_0 u_t)_\cO \\
&= & \|u_t\|^2_\Om  - \|A^{1/2}u\|^2_\Om
-(f(u),u)_\Om -(\ga_\Om v, u)_\Om \\
&&
+(\nu \Delta v+G, N_0 u)_\cO+(v, N_0 u_t)_\cO\\
&= & \|u_t\|^2_\Om  - \|A^{1/2}u\|^2_\Om
-(f(u),u)_\Om  \\
&&
-\nu (\nabla v, \nabla N_0 u)_\cO+ (G, N_0 u)_\cO+(v, N_0 u_t)_\cO.
\end{eqnarray*}
Here above we also use the fact
that $A_0(v-N_0u_t)=-P_L\Delta v$ and the Green formula:
\[
\nu(\Delta v, N_0 u)_\cO=-\nu (\nabla v, \nabla N_0 u)_\cO+(\ga_\Om v, u)_\Om.
\]
Therefore using Proposition~\ref{pr:stokes} we obtain that
\begin{eqnarray}\label{lyap-d}
& &\frac{d}{dt} \left[ (u,u_t)_\Om +(v, N_0 u)_\cO\right] \\
&\le & 2\|u_t\|^2_\Om  - (1- \delta) \|A^{1/2}u\|^2_\Om
-(f(u),u)_\Om
+c_\delta \left(\|\nabla v\|^2_\cO+ \|G \|^2_\cO\right)
\nonumber
\end{eqnarray}
for every $\delta>0$.
By our hypotheses (see \eqref{phi-2}) there is $c^*_f>0$ and $c_f\ge0$ such that
\begin{equation}\label{phi-3}
(f(u),u)_\Om - c_f^* \int_\Om \Phi(u) d\Om+\frac12\|A^{1/2}u\|^2_\Om\ge - c_f.
\end{equation}
Therefore from \eqref{lyap-d} we have that
\begin{eqnarray*}
& &\frac{d}{dt} \left[ (u,u_t)_\Om +(v, N_0 u)_\cO\right] \\
&\le & 2\|u_t\|^2_\Om  - \frac14 \|A^{1/2}u\|^2_\Om
-c_f^* \int_\Om \Phi(u) d\Om
+c_{1/4} \left(\|\nabla v\|^2_\cO+ \|G \|^2_\cO\right)+c_f.
\end{eqnarray*}
After selecting $\eta>0$ small enough,
this  relation along with \eqref{energy-d1} implies that
\begin{equation}\label{energy-wd1}
\frac{d}{dt}\sW(v(t), u(t), u_t(t))+c_0\sW(v(t), u(t), u_t(t))\le
c_1\left( c_f+  \|G\|_\cO^2\right)
\end{equation}
with some $c_0,c_1>0$.
This yields
\begin{equation}\label{energy-wd1a}
\sW(v(t), u(t), u_t(t))\le \sW(v_0, u_0, u_1)e^{-c_0t}
+\frac{c_1}{c_0}\left(c_f+ \|G\|_\cO^2\right)
\left[ 1-e^{-c_0t}\right].
\end{equation}
For $\eta>0$ small enough we definitely have that
\begin{equation}\label{energy-wd3}
\frac12\cE(v, u, u_t)\le\sW(v, u, u_t)\le 2 \cE(v, u, u_t).
\end{equation}
Therefore the standard argument
(see, e.g., \cite{BabinVishik,Chueshov,Temam})
 yields dissipativity with a forward invariant bounded absorbing set.
\end{proof}
The argument given above allows us to obtain the following assertion
on asymptotic stability.
\begin{proposition}\label{pr:as-stab}
If $G\equiv 0$ and relation \eqref{phi-2} holds with $c_2(\delta)\equiv0$,
Then there exist $c_R>0$ and  $\al>0$ such that
\begin{equation}\label{exp-stab-1}
\|S_tU\|_\cH\le c_R e^{-\al t}~~\mbox{for any $U\in\cH$ such that $\|U\|_\cH\le R$}.
\end{equation}
This means that in this case the global attractor $\Ac$ consists
of a single point, $\Ac=\{(0;0;0)\}$, which is exponentially attractive.
\end{proposition}
\begin{proof}
In the case considered  we have $c_f=0$ in \eqref{phi-3} and also  $G\equiv0$.
 Thus  \eqref{energy-wd1a} has the form
\[
\sW(v(t), u(t), u_t(t))\le \sW(v_0, u_0, u_1)e^{-c_0t},~~ t>0.
\]
Hence \eqref{energy-wd3} implies \eqref{exp-stab-1}.
\end{proof}
We note that  the hypotheses of Proposition~\ref{pr:as-stab} holds true
if $G\equiv 0$ and relation \eqref{phi-exm}
is valid with $\psi_i\in C^2(\R)$
possessing the additional properties: $\psi_i(s)\ge \psi_i(0)=0$ and $\psi_i'(s)$ is non-decreasing
function.
\par
As a consequence of Proposition~\ref{pr:as-stab}
for $f(u)\equiv 0$ we have
 the following assertion, which enforces Proposition~\ref{pr:cA-generator}.
\begin{corollary}\label{co:as-stab}
The operator $\cA$ given by \eqref{opA-1} and \eqref{opA-2}
generates an exponentially stable  $C_0$-semigroup  of contractions $e^{-t\cA}$ in $\cH$.
In particular, there exists positive $C,\al$ such that
\begin{equation}\label{exp-stab}
\|e^{-t\cA} U\|_{L(\cH,\cH)}\le C e^{-\al t}~~\mbox{for all}~t>0.
\end{equation}
\end{corollary}
This  corollary improves the result in \cite{Ggob-mmas09} which states the strong stability
only.
\subsection{Quasi-stability}
We use the method  developed in  \cite{cl-jdde}  (see also \cite{cl-mem,cl-book}
and the references therein) to obtain asymptotic smoothness and
finiteness of fractal dimension of the attractor.
This method based on quasi-stability properties of the system and
involves
the so-called stabilizability estimate.
\begin{proposition}[Stabilizability estimate]\label{pr:st-est}
Let
\[
S_tU_0\equiv U(t)=( v(t); u(t); u_t(t))~~ and
 ~~ S_tU_0^*\equiv U^*(t)=( v^*(t); u^*(t); u^*_t(t))
\]
be two semi-trajectories such $\|S_tU_0\|_\cH,\|S_tU_0^*\|_\cH\le R$
for all $t>0$ and for some $R>0$.
Then there exist  positive constants $c_0$, $\om$ and $c_R$ such that
\begin{equation}\label{stab-est}
\|S_tU_0-S_tU_0^*\|_\cH\le c_0 e^{-\om t} \|U_0-U_0^*\|_\cH
+c_R\int_0^t e^{-\om( t-\tau)} \|u(\tau)-u^*(\tau)\|_\Om d\tau
\end{equation}
for all $t>0$.
\end{proposition}
\begin{proof}
Using the representation \eqref{mild-s}, the exponential stability
of the  linear  semigroup $e^{-t\cA}$  given by \eqref{exp-stab},
 and also the local Lipschitz property in \eqref{Lip-f}
 we obtain that
\begin{equation}\label{stab-est1}
\|S_tU_0-S_tU_0^*\|_\cH\le C e^{-\al t} \|U_0-U_0^*\|_\cH
+C_R\int_0^t e^{-\al( t-\tau)} \|z(\tau)\|_{\si,\Om} d\tau
\end{equation}
for all $t>0$, where $z(t)=u(t)-u^*(t)$ and $0<\si<1$, By interpolation,
\[
\|z\|_{\si,\Om}\le \delta \|A^{1/2}z\|_{\Om}+ C_\delta\|z\|_{\Om}
\le \delta \|S_tU_0-S_tU_0^*\|_\cH+ C_\delta\|z\|_{\Om}
\]
for every $\delta>0$. Substituting this relation in \eqref{stab-est1}
and applying the Gronwall type argument with an appropriate choice $\delta$
we obtain  \eqref{stab-est}.
\end{proof}
We note that the property stated in Proposition \ref{pr:st-est}
means that the system $(\cH,S_t)$ is quasi-stable (in the sense
of \cite[Definition 7.9.2]{cl-book}). This observation make it possible
to obtain several important dynamical properties
at the abstract level.

\subsection{Completion of Theorem~\ref{th:main}}
Using the stabilizability estimate in \eqref{stab-est}
via the Ceron--Lopes type criteria (see
\cite[Corollary~2.7]{cl-mem}) we can easily prove
that $(\cH,S_t)$ is asymptotically smooth
(for some details for similar systems see, e.g., \cite{cl-mem,cl-amo08}
and also \cite[Section 7.9]{cl-book} at the abstract level).
Thus there exists a compact global attractor $\Ac$ for $(\cH,S_t)$.
\par
This attractor has a finite fractal dimension. To see this
one should apply the same argument as in \cite[Theorem 4.3]{cl-mem},
see also the argument given in the proof of Theorem 4.1
in \cite{cl-amo08} for the case of  thermoelastic plate models.
\par
The regularity property in \eqref{reg}
follows by the same argument as in the proof of \cite[Theorem 4.17]{cl-mem},
see also the proof of Theorem 6.2 \cite{cl-amo08}
and argument given in \cite[Section 7.9.2]{cl-book}
for the abstract quasi-stable systems.
\par
Thus the proof of Theorem~\ref{th:main} is complete.

\begin{remark}
{\rm
As in \cite{cl-mem} and \cite[Section 7.9]{cl-book}
we can use the stablizability estimate to construct
exponential fractal attractor (whose dimension is finite
in some extended space) and
prove the existence of finite number of determining functionals
supported by displacement component of the plate.
We also note that using the same approach as in \cite[Theorem 4.23]{cl-mem}
(see also \cite[Section 8.7]{cl-book}) and assuming additional smoothness
of the potential $\Phi$
 we can obtain higher regularity
of time derivatives of the trajectories  from the attractor.
}
\end{remark}

\end{document}